\documentclass[11pt]{article}
\usepackage{amsmath}
\usepackage{amsthm}
\usepackage{mathtools}
\usepackage{amssymb}

\usepackage{enumitem}
\usepackage{hyperref}
\hypersetup{pdfcreator=XeLaTeX with hyperref}

\usepackage{geometry}
\usepackage{setspace}
\title{\textbf{{\normalsize Approximation on compact sets of functions and all derivatives}}} 
\author{\textsc{Sotiris Armeniakos, Giorgos Kotsovolis and Vassili Nestoridis} \\ \\ \textit{National and Kapodestrian University of Athens} \\ \textit{Department of Mathematics}}
\date{}

\begin{document}
\maketitle 

\begin{abstract}
\noindent
In Mergelyan type approximation we uniformly approximate functions on compact sets $K$ by polynomials or rational functions or holomorphic functions on varying open sets containing $K$. In the present paper we consider analogous approximation, where uniform convergence on $K$ is replaced by uniform approximation on $K$ of all order derivatives.
    
\end{abstract}
\noindent
AMS classification numbers: $30E10,30H50$ \\ \\
Key words and phrases: Mergelyan theorem, uniform approximation of all order derivatives, $C^\infty \, ,A^\infty \, ,A(K) \, ,C(K)$, holomorphic functions, completion of a metric space.

\section{Introduction}
Let $K$ be a compact set in the complex plane $\mathbb C$, or more generally in $\mathbb C^d$. By $O(K)$ we denote the set of all complex functions $f$, such that there exists an open set $V_f$ containing $K$ with the property that $f$ is defined and is holomorphic on $V_f$. One endows $O(K)$ with the supremum norm on  $K$. Then the completion of the metric space is denoted by $\overline{O}(K)$ and coincides with the closure of $O(K)$ in $C(K)$. It is well known hat $\overline{O}(K)$ is contained in the classical algebra $A(K)=\{f:K\longrightarrow \mathbb C$ continuous on $K$ and holomorphic in $K^o$ \}. If $K^o=\varnothing$, then $A(K)=C(K)$. Mergelyan type theorems are those where under some assumptions on $K$ we can conclude that $A(K)=\overline{O}(K)$. In particular, if $K\subset \mathbb{C}$ and $\mathbb C-K$ is connected, then every function f in $A(K)$ can be uniformly on $K$ approximated by polynomials ;thus, $A(K)=\overline{O}(K)$ in this case. More generally, if $K\subset \mathbb C$ and $\mathbb C-K$ has a finite number of components, then every $f$ in $A(K)$ can be uniformly on $K$ approximated by rational functions with prescribed poles, one in each component of $\mathbb C_\infty - K$; thus, in this case we also have $A(K)=\overline{O}(K)$. \\
For $K\subset \mathbb C$ a characterization of whether it holds that $A(K)=\overline{O}(K)$ or not has been obtained by Vitushkin using continuous analytic capacity. (\cite{viruvskin1967analytic}, \cite{vitushkin1975uniform}) \\
In the present paper, for $K$ a compact planar set, we replace uniform on $K$ approximation by uniform on $K$ approximation of all order derivatives. More precisely, we endow $O(K)$ with the metric $$d(f,g)= \sum\limits_{k=0}^{\infty}\dfrac{1}{2^{k}}\dfrac{\sup\limits_{z\in K}\vert f^{(k)}(z) - g^{(k)}(z)  \vert}{1+\sup\limits_{z\in K}\vert f^{(k)}(z) - g^{(k)}(z)\vert}$$ ,where $F^l$ denoted the $l^{th}$ derivative of $F\in O(K)$, which is well defined on $K$ because $F$ is holomorphic on an open set containing $K$. Now the analogue of $\overline{O}(K)$ is the completion $B(K)$ of the metric space $(O(K),d)$. Every element of $B(K)$ is a sequence $g=(g_l)_{l=0}^\infty \in A(K)^{\aleph o}$. Our effort is to determine explicitly which subset of $A(K)^{\aleph o}$ is $B(K)$. At the end of the paper we show that, if $K$ is denumerable then $B(K)=A(K)^{\aleph o}$. If $K=\{z\in C: \vert z \vert \leq1\}$, then $B(K)= \{g=(g_l)_{l=0}^\infty:$ there exists a holomprhic function $f:\{z\in \mathbb C: \vert z \vert <1\} \longrightarrow \mathbb C$, such that each derivative $f^l$ extends continuously on $\{z\in \mathbb C: \vert z \vert \leq1\}$ and $g_l=f^l$ \}. Roughly speaking $B(K)$ can be identified with $A^\infty(D)$, $D$ being the unit disk. If $K=[0,1]$ then $B(K)$ can be identified with $C^\infty[0,1]$. \\
In the general case we define a subset $\Gamma(K)$ of $C(K)^{\aleph o}$ containing $B(K)$ and we give sufficient conditions on $K$ so that $B(K)=\Gamma(K)$. This happens in all previously mentioned cases. We also prove a lemma stating that rational functions with poles outside $K$ are always ''dense'' in $B(K)$, where $f=(f^0,f^1,f^2,...)$. \\
Extensions in several variables are possible and will be treated in future papers. When $K\subset \mathbb C^d,d>1$, then $A(K)$ will be replaced by a smaller algebra $A_D(K)=\{f\in C(K):f$ is holomorhpic on any analytic disk contained in $K$, even meeting the boundary $K\}$, which satisfies $\overline{O}(K)\subset A_D(K)\subset A(K)$. (\cite{falco2019function}). The algebra $A_D(K)$ is more appropriate for approximation in several complex variables than the classical algebra $A(K)$.

\section{A function Algebra and an Important Lemma}
In this section we introduce an algebra of functions, which will be of interest throughout this paper. Namely, we have the standard definition: \\

\noindent
\textbf{Definition 2.1} Let $K \subsetneq \mathbb C$ be a planar compact set. We call $O(K)$ the set of all functions $f$ such that there can be found an open set $V_{f} \supseteq K$ (dependent on $f$) , such that $f$ is holomorphic on $V_{f}$.\\ \\ 
Typically this set is endowed with the supremum norm on $K$ and is thus considered as a subset of $C(K)$, the space of continuous functions on $K$. However, in this paper we will provide this space with another topology. We give the next definition: \\ \\
\noindent
\textbf{Definition 2.2} We endow the space $O(K)$ of Definition 1, with the metric introduced by the denumerable family of seminorms $\sup\limits_{z\in K} \vert {f^{i}(z)} \vert$ where $f^{i}$ denotes the derivative of i-order of $f$. Namely $O(K)$ is given the metric $$d(f,g)= \sum\limits_{k=0}^{\infty}\dfrac{1}{2^{k}}\dfrac{\sup\limits_{z\in K}\vert f^{(k)}(z) - g^{(k)}(z)  \vert}{1+\sup\limits_{z\in K}\vert f^{(k)}(z) - g^{(k)}(z)\vert}$$ .\\ 
From now on we will refer to this particular metric as $d$. 

\noindent
Notice that Definition 2.2 is well posed, since each function in $O(K)$ is defined and holomorphic in a neighborhood of $K$. \\ \\
\noindent
\textbf{Definition 2.3} We call $B(K)$ the completion of $O(K)$ with respect to the metric introduced in Definition 2.2 .\\ \\
\noindent
\textbf{Remark 2.4} Notice  that taking an element $f$ of $O(K)$, for some compact set $K\subset \mathbb C$, we can define the sequence $(f,f',f'',...)$ of its derivatives and call this sequence, the sequence corresponding to $f$. Note now that $d(f_n,f)\rightarrow 0$ for some functions $f_n,f$ of $O(K)$ if and only if the corresponding sequence of $f_n$ converges on the corresponding sequence of $f$, where the limit is realized as the term to term limit via the supremum norm on $K$. These sequences, which correspond to elements of $O(K)$ are of course sequences in $A(K)^{N_o}$. In this way we can refer to elements of the complete space $B(K)$ as sequences in $A(K)^{N_o}$. Often we identify a function $f$ with the sequence of its derivatives. Thus, we can say that a subset of $O(K)$ is also a subset of $B(K)$. \\ \\
\noindent
In this paper one of the questions we are interested is in what length can we determine the set $B(K)$ explicitly. To answer this we will use a lemma, making use of Runge's Theorem, which we now state. \\ \\
\noindent
\textbf{Theorem (Runge)} \textit{Let $K \subsetneq \mathbb C$ be a compact planar set and $L$ be a set, containing at least one point for every connected component of the set $\mathbb C_\infty-K$. Then every function $f$ in $O(K)$ can be uniformly on $K$ approximated by rational functions with poles in $L$.} \\ \\
We also give another version of Runge's Theorem: \\ \\
\textbf{Theorem (Runge)} \textit{Let $U \subsetneq \mathbb C$ be an open planar set and $L$ be a set, containing at least one point from every connected component of the set $\mathbb C_\infty-U$. Then every function $f$ in $H(U)$ can be approximated uniformly on compact subsets of $U$ by rational functions with poles in $L$.} \\ \\
Runge's Theorem gives us that given a compact set $K$ and a set $L$ as that in the theorem, $R_L$ is dense in $O(K)$ with respect to the supremum norm on $K$,where $R_L$ denotes the set of rational functions with poles only in $L$. However, we have already stated that throughout this paper the metric on $O(K)$ that concerns us is not the supremum one, but the metric $d$. It is therefore reasonable to ask whether $R_L$ is dense in $O(K)$ with respect to d. The answer to that is affirmative and for the proof of this we will need the following lemma: \\ \\ 
\noindent
\textbf{Lemma 2.5} (\cite{diamantopoulos2006universal})  \textit{Let $K$ be a compact planar set and $L$ a set containing one point from every connected component of $C_\infty-K$. If V is an open set containing $K$, there exists an open set $W$, with $K\subset W\subset V$, such that every connected component of $\mathbb C_\infty-W$ contains a point of $L$.} \\ \\
We now prove our proposition: \\ \\
\textbf{Proposition 2.6} \textit{Let $K\subset \mathbb C$ be a compact planar set and $L$ be a set, containing at least one point from every connected component of $\mathbb C_\infty-K$. Let $R_{L}$ be the set of rational functions with poles in $L$. $(R_{L} \subset O(K))$. Then $R_{L}$ is dense in $B(K)$.} \\ \\ 
\noindent
\textit{Proof.} It is, of course, sufficient to show that $R_L$ is dense in $O(K)$. Let $f \in O(K)$. Thus, there exists an open set $V_f \supset K$ such that $f$ is analytic on $V_f$. By Lemma 2.5 we can find an intermediate open set $W$, with $K\subset W\subset V_f$, such that every connected component of $\mathbb C_\infty-W$ contains a point of $L$. By Runge's Theorem (second version) f can be uniformly approximated on compact subsets of $W$ by rational functions $r_n$ with poles in $L$. By Weierstrass' Theorem $f^l$ is approximated by $r_n^l$ uniformly on $K$ for every $l\in N$. Hence, we have our result. \qedsymbol \\ \\ 

\noindent
We will also state a form of Mergelyan's Theorem for future reference. \\

\noindent
\textbf{Theorem (Mergelyan) } (\cite{rudin2006real})  \textit{Let $K$ be a compact planar set such that $\mathbb C_\infty-K$ has finitely many connected components. Let $L$ denote a set containing one point from each connected component of $\mathbb C_\infty-K$. Then the rational functions with poles in $L$ are uniformly dense in $A(K)=\{f:K\longrightarrow \mathbb C$ continuous on $K$ and holomorphic in $K^o\}$. }

\section{The function algebra $\Gamma(K)$ and its relation with $B(K)$}

In this section we will introduce a new function algebra, which has the advantage of being intuitively clearer than $B(K)$. To motivate for the following definition we noted in Remark 2.4 that taking an element $f$ of $O(K)$, for some compact set $K\subset \mathbb C$, we can define the sequence $(f,f',f'',...)$ of its derivatives and call this sequence, the sequence corresponding to $f$. These sequences, which correspond to elements of $O(K)$ are sequences in $A(K)^{N_o}$. A property of them is that taking $f\in O(K)$ and $\gamma$ a rectifiable curve in $K$ starting at a point $a$ and ending at a point $b$ and $\phi$ a holomorphic function in some neighborhood of $\gamma$ , then  
\begin{equation}
\int_{\gamma}\phi(z){f}^{l+1}(z)\,dz\,=\,f^{l}(z)\phi (z)\vert{^b _a}-\int_{\gamma}f^{l}(z)\phi'(z)\,dz
\end{equation}
Now if we take $g \in B(K)$, it can be expressed as a sequence $(g_0,g_1,...) \in A(K)^{N_0}$, where there exist $f_n \in O(K)$ such that for every $k \in N$

\[\sup\limits_{z \in K} \vert {f_n}^{k}(z) - g_k(z)\vert \rightarrow 0\]
By (1):
$$\int_{\gamma}\phi(z){f_n}^{l+1}(z)\,dz\,=\,f_{n}^{l}(b)\phi (b)-f_{n}^{l}(a)\phi (a)-\int_{\gamma}f_{n}^{l}(z)\phi'(z)\,dz$$
Taking $n\rightarrow\infty$ we obtain:
$$\int_{\gamma}\phi(z){g_{l+1}}(z)\,dz\,=\,g_{l}(b)\phi (b)-g_{l}(a)\phi (a)-\int_{\gamma}g_{l}(z)\phi'(z)\,dz$$
Keeping that in mind we define: \\ \\
\textbf{Definition 3.1} Let $K$ be a planar compact set. We define $\Gamma(K)$ as the set of sequences $(g_k)_{k=0}^\infty$ in $A(K)^{\aleph_0}$ such that for every rectifiable curve $\gamma$ in $K$ starting at some point $a$ and ending at some point $b$, and every $\phi$ holomorphic in some neightbboorhood of $\gamma$ we have that 
$$\int_{\gamma}\phi(z){g_{l+1}}(z)\,dz\,=\,g_{l}(b)\phi (b)-g_{l}(a)\phi (a)-\int_{\gamma}g_{l}(z)\phi'(z)\,dz$$ for every $k \in N$. \\ \\
 \textbf{Remark 3.2} We provide $\Gamma(K)$ with the metric $d$, where $$d((f_k)_{k=0}^\infty,(g_k)_{k=0}^\infty)= \sum\limits_{k=0}^{\infty}\dfrac{1}{2^{k}}\dfrac{\sup\limits_{z\in K}\vert f_k(z) - g_k(z)  \vert}{1+\sup\limits_{z\in K}\vert f^{(k)}(z) - g^{(k)}(z)\vert}$$
 Due to the conversation preceeding Definition 3.1, we have that $B(K) \subset\Gamma(K)$ and the inclusion map is an isometry. \\
 We will now show that in many familiar cases $B(K)=\Gamma(K)$. \\ \\ 
 \textbf{Theorem 3.3} \textit{Let $K\subset \mathbb C$ be a compact set, such that: \\
 \begin{enumerate}
     \item $\mathbb C_\infty-K$ has a finite number of connected components, $V_0, V_1,...,V_n$ and let $a_0=\infty \in V_0$ and $a_j \in V_j , \, j=0,1,2..,n$ be fixed.
     \item There exists $S \subset K$ with $\overline{S}=K, a\in S$ and $M<\infty$, so that for all $z\in S \,\,$ there exists a rectifiable curve $\gamma _{a,z}$ in $K$ starting at $a$ and ending at $z$ with length $l(\gamma_{a,z})\leq M$  
     \item There exist $\delta_i $ in $ K, \, i=1,2,...,n$ closed rectifiable curves such that $Ind(\delta_i,a_j)=0$ for $j\neq i$ and $Ind(\delta_i,a_i)=1$.
     \end{enumerate}
     Then $\Gamma(K)=B(K)=\overline{{R_L}}$ where $L=\{a_0,a_1,...a_n\}$ and the closure of $R_L$ is taken in $B(K)$.}\\ \\
\textbf{Proof.} Let $g=(g_0,g_1,g_2,..) \in \Gamma(K)\,$ and $N\in \mathbb{N}$. It suffices to find $f_n\in R_L$ such that \\ $\sup\limits_{z\in K}\vert f_n^l(z)-g_l(z)\vert\xrightarrow[n\rightarrow \infty]{}
  0$ for $l=0,1,2...N$. By Mergelyan's Theorem $\exists\, h_m \, \in R_L$, such that \\ $\sup\limits_{z \in K}\vert h_m(z)-g_N(z)\vert \xrightarrow[m\rightarrow \infty]{} 0$.
Now define $h_{m,r,i}$ for $1\leq r\leq N$ and $1\leq i\leq n$ to be $h_{m,r,i}=\underset{z=a_i}{Res}\,h_m(z){(z-a_i)}^{r-1}$.
Take $\widetilde{h}_m(z)=h_m(z)-\sum\limits_{r=1}^{N}\sum\limits_{i=1}^{n} \dfrac{h_{m,r,i}}{(z-a_i)^r}$. $\widetilde{h}_m \in R_L$ easily.\\  
By construction $\widetilde{h}_m$ has now an $N^{th}$ order primitive, that is there exists $H_m$ in $R_L$
 such that $H_m^N(z)=\widetilde{h}_m(z)$. We can also take $H_m^r(a)=g_r(a)$ for $0\leq r \leq N-1$. We will now show that $$\sup\limits_{z \in K} \vert H_m^{(k)}(z) - g_k(z)\vert \xrightarrow[m\rightarrow \infty]{} 0$$
 for $k=0,1,2,...,N$.\\
 We will first show that $\big\vert \dfrac{h_{m,r,i}}{(z-a_i)^r}\big\vert \rightarrow 0$ uniformly on $K$, or equivalently, since $dist(a_i,K)>0$, that $h_{m,r,i} \rightarrow 0$. To see that notice that: $$2\pi i\underset{z=a_i}{Res}h_m(z)(z-a_i)^{r-1}=\int_{\delta_i}h_m(z)(z-a_i)^{r-1}\,dz \rightarrow\int_{\delta_i}g_N(z)(z-a_i)^{r-1}\,dz=-(r-1)\int_{\delta_i}g_{N-1}(z)(z-a_i)^{r-2}\,dz$$ as $m\rightarrow \infty$, since $\delta_i$ is closed and $g\in \Gamma(K)$. By induction: $$\int_{\delta_i}g_N(z)(z-a_i)^{r-1}=(-1)^{(r-1)}(r-1)!\int_{\delta_i}g_{N-r}(z)\,dz=0$$ for $0\leq r\leq N-1$. Therefore, $h_{m,r,i} \rightarrow 0$ as $m\rightarrow \infty$. \\
 Now $\sup\limits_{z \in K} \vert h_m(z)-g_N(z) \vert \xrightarrow[m\rightarrow \infty]{} 0$ and $\sup\limits_{z \in K} \vert \widetilde{h}_m(z)-g_N(z) \vert \xrightarrow[m\rightarrow \infty]{}0$. \\
 Therefore we have that $\sup\limits_{z \in K} \vert H_m^N(z)-g_N(z) \vert \xrightarrow[m\rightarrow \infty]{} 0$. We proceed by induction: \\
 Suppose for $ r\in {1,2,...,N-1}$ that $\sup\limits_{z \in K} \vert H_m^{r+1}(z)-g_{r+1}(z) \vert \rightarrow 0$.
 Then we have that taking $b \in S$ and $\gamma_{a,b}$ in $K$ such that $l(\gamma_{a,b})\leq M$:\\
 $$\int_{\gamma_{a,b}}H_m^{r+1}(z)\,dz=H_m^{r}(b)-H_m^{r}(a)$$ and 
 $$\int_{\gamma_{a,b}}g_{r+1}(z)\,dz=g_{r}(b)-g_{r}(a)$$. Thus, since $H_m^r(a)=g_r(a)$ we have:
 $$\vert H_m^r(b) - g_r(b)\vert \leq \int_{\gamma_{a,z}}\vert H_m^{r+1}(z) - g_{r+1}(z)\vert \,\vert dz \vert\leq M\cdot \sup\limits_{z \in K}\vert H_m^{r+1}(z) - g_{r+1}(z)\vert \rightarrow 0$$, as $m\rightarrow \infty$, for all $b\in S$ and $0\leq r\leq N-1$.\\ Since $\overline{S}=K$, we have the proof completed. $\qedsymbol$ \\ \\
 We now make a standard definition: \\ \\
 \textbf{Definition 3.4} Let $\Omega \in \mathbb C$ be a bounded domain such that ${\overline{\Omega}^o=\Omega}$. A holomorphic function $f:\Omega \longrightarrow \mathbb C$ belongs to $A^{\infty}(\Omega)$ if for every $l=0,1,2...$
 the derivative $f^l$ can be continuously extended on $\overline{\Omega}$. The topology in $A^{\infty}(\Omega)$ is induced by the seminorms $\sup\limits_{z \in \overline{\Omega}}\vert f^l(z)\vert$. \\ \\
 As Corollaries of Theorem 3.3, we have the following statements: \\ \\
 \textbf{Corollary 3.5} \textit{Let $\Omega \subset \mathbb C$ be a bounded domain such that ${\overline{\Omega}^o}=\Omega$ and such that $\mathbb C_\infty-\overline{\Omega}$ is connected. We also assume that there exists $M<\infty$, such that for every $a,b \in \Omega$ there exists a rectifiable curve $\gamma_{a,b}$ in $\Omega$ joining $a$ and $b$ with length $l(\gamma_{a,b})\leq M$. Then, if $P$ denotes the set of polynomials we have that $\overline{P}=A^{\infty}(\Omega)=B(\overline{\Omega})=\Gamma(\overline{\Omega})$, where the closure of $P$ is taken in the topology of $A^\infty(\Omega)$.}\\
 \textbf{Proof.} Setting $K=\overline{\Omega}, S=\Omega$ in Theorem 3.3, we have that $\overline{P}=B(\overline{\Omega})=\Gamma(\overline{\Omega})$. Since $P \subset A^{\infty}(\Omega) \subset \Gamma(\overline{\Omega})$ and the inclusion map $A^{\infty}(\Omega) \subset \Gamma(\overline{\Omega}) $ is an isometry, the statement is proved. $\qedsymbol$\\ \\
 \textbf{Corrolary 3.6} \textit{Let $\Omega \subset \mathbb C$ be a bounded domain such that ${\overline{\Omega}^o=\Omega}$ and such that $\Omega$ is bounded by a finite set of disjoint Jordan curves. Let $L=\{a_0,a_1,a_2,...a_n\}$ containing exactly one point from each connected component of $\mathbb C_\infty-\overline{\Omega}$ with $a_0=\infty$.We also assume that there exists $M<\infty$, such that for every $a,b \in \Omega$ there exists a rectifiable curve $\gamma_{a,b}$ in $ \Omega$ joining $a$ and $b$ with length $l(\gamma_{a,b})\leq M$. Then $\overline{R_L}=A^{\infty}(\Omega)=B(\overline{\Omega})=\Gamma(\overline{\Omega})$, where the closure of $R_L$ is in the topology of $A^\infty(\Omega)$.}\\ \\
 \textbf{Proof.} The proof is identical to that of Corollary 3.5 . $\qedsymbol$ \\ \\
 \textbf{Remark 3.7} If $K$ is a compact rectifiable curve $\gamma:[0,1]\rightarrow \mathbb C$, such that $\mathbb C_{\infty}-K$ has a finite set of connected components, then the assumptions of Theorem 3.3 are satisfied. Therefore, if we pick a set $L$, containing a point from each connected component of $\mathbb C_{\infty}-K$, then $\overline{R_L}=B(K)=\Gamma(K)$. Now, suppose in addition the curve K is locally one to one and satisfies the following: \\
 For every $a \in [0,1]$ there exists a constant $C_a$ denendent on $a$ and a $\delta>0$, such that, for every $t\in(a-\delta,a+\delta)$ we have that: $$l(a,t)\leq C_a\dot\vert \gamma(t)-\gamma(a)\vert$$,where $l(a,t)$ is the length of $\gamma_{\vert[a.t]}$ if $a<t$ or of $\gamma_{\vert[t.a]}$ if $t\leq a$.\\ Then condition (1) implies that for a sequence $g_l$ in $B(K)$ we have : $$\int_a^t{g_{l+1}}(z)\,d\gamma \,=\,g_{l}(\gamma(t))-g_{l}(\gamma(a))$$ Thus taking derivatives with respect to the position we have:
 $$\begin{aligned}\frac{dg_l}{dz}\vert_{z=\gamma(a)}=\lim_{t\rightarrow a}\frac{g_{l}(\gamma(t))-g_{l}(\gamma(a))}{\gamma(t)-\gamma(a)}=\lim_{t\rightarrow a}\frac{\int_a^t{g_{l+1}}(z)d\gamma}{\gamma(t)-\gamma(a)}=\lim_{t\rightarrow a}\frac{\int_a^t{g_{l+1}}(z)-g_{l+1}(\gamma(a))d\gamma}{\gamma(t)-\gamma(a)}+g_{l+1}(\gamma(a))=g_{l+1}(\gamma(a))\end{aligned}$$ since $$\limsup\limits_{t\rightarrow a}\vert\frac{\int_a^t{g_{l+1}}(z)-g_{l+1}(\gamma(a))\,d\gamma}{\gamma(t)-\gamma(a)}\vert \leq \limsup\limits_{t\rightarrow a}\frac{l(a,t)\cdot \sup\limits_{\vert s-a\vert\leq \vert t-a\vert}\vert g_{l+1}(\gamma(s)) - g_{l+1}(\gamma(a))\vert}{\vert \gamma(t)-\gamma(a)\vert } \leq $$ $$\leq C_a\cdot \limsup\limits_{t\rightarrow a}\sup\limits_{\vert s-a\vert\leq \vert t-a\vert}\vert g_{l+1}(\gamma(s)) - g_{l+1}(\gamma(a))\vert\rightarrow 0$$ due to the continuity of $g_{l+1}$.\\
 We thus get that the the differentiation of $g_l$ with respect to the position yields $g_{l+1}$. It follows that $\overline{R_L}=B(K)=\Gamma(K)$ coincides with the set $C^{\infty}(K)$ of all functions $f:K\longrightarrow \mathbb C$, which have derivatives of any order on $K$ with respect to the position.\\ \\
 \textbf{Remark 3.8} Using Remark 3.7, we see that if $K=[a,b]$, then $B(K)=C^{\infty}[a,b]$ and if $K=\{e^{i\theta}:\theta\in R\}=T$, then $B(T)=C^{\infty}(T)$.\\
 \textbf{Remark 3.9} Using the previous results we can make deductions about $B(K)$ on more complex situations. For example take $K=K_1\cup K_2 \cup K_3$, where $K_1$ and $K_3$ are disjoint closed disks and $K_2=[A,B]$, where $A\in \partial K_1, \, B\in \partial K_3$ and $K_2$ is disjoint from $K_1$ and $K_3$, except for the points $A$ and $B$. Then $f:K \longrightarrow \mathbb C$ belongs to $B(K)$ if and only if $f_1=f_{\vert K_1} \in B(K_1), f_2=f_{\vert K_2} \in B(K_2),f_3=f_{\vert K_3} \in B(K_3)$ and for every $l=0,1,2,3,... \, f_1^l(A)=f_2^l(A)$ and $ f_3^l(B)=f_2^l(B)$. We know, however, that $B(K_1)=A^{\infty}(K_1),B(K_3)=A^{\infty}(K_3)$ and $B(K_2)=C^{\infty}(K_2)$. We thus have an exact description of $B(K)$. \\ \\
 We now give a result, concerning the extendability of functions $g_k, k=0,1,2,3,..$ where $(g_k)\in B(K)$, $K$ being a compact set with no isolated points. \\ \\
 \textbf{Theorem 3.10} \textit{Let $K \subset \mathbb C$ be compact without isolated points. Then there exists $(g_k)_{k=0}^{\infty} \in B(K)$, such that for every $k=0,1,2,..$ and every open disk $D$, such that $D\cap K\neq\varnothing, D\cap K^c\neq\varnothing$, there is no holomorhpic  function $F:D \longrightarrow \mathbb C$ such that $F_{\vert D\cap K}={g_k}_{\vert D\cap K}$. The set of such sequences $(g_k)_{k=0}^\infty$ is $G_{\delta}$-dense in $B(K)$.} \\
 \textbf{Proof.} Define $E_k=\{(g_n)_{n=0}^{\infty}\in B(K)$, such that for every open disk $D$ such that $D\cap K\neq\varnothing,D\cap K^c\neq\varnothing$ there is no holomorphic function $F:D\longrightarrow \mathbb C$ 
, such that  $F_{\vert D\cap K}={g_k}_{\vert D\cap K}\}$. We want to show that $\bigcap\limits_{k=0}^{\infty}E_k$ is $G_\delta$-dense. Since $B(K)$ is a complete metric space, by Baire's Theorem we need to show $E_k$ is $G_\delta$-dense for each $k=\{0,1,2,..\}$. \\
Let $k\in N$. Notice that for $(g_n)_{n=0}^\infty \in (E_k)^c$, there exists $D$, open disk such that $g_k$ is holomorphically extendable on $D$ and $D\cap K\neq \varnothing,D\cap K^c\neq \varnothing$. If we pick a smaller disk $D' \subsetneq D$, such that $D'\cap K\neq \varnothing,\,D'\cap K^c\neq \varnothing$, then the extension on $D'$ is also bounded. \\
Thus $(E_k)^c=\bigcup\limits_{n\in N}^\infty \{ (g_s)_{s=0}^\infty \in B(K)$, such that there exists an open disk $D$, with $D\cup K\neq \varnothing,\,D\cap K^c\neq \varnothing$ and a holomorphic function $F$ on $D$, such that $F_{\vert D\cap K}={g_k}_{\vert D\cap K}$ and $\vert F(z)\vert\leq n$ for $z \in D$ for some natural number $n$\}. It is known that disks with rational centers and rational ratios form a subbase of the topology of $\mathbb C$. \\ Assume that $r_1,r_2,...$ is an enumeration of the rational numbers and define $S=\{(k_1,k_2,k_3)$, such that $r_{k_3}>0$ and $D\cap K\neq \varnothing, D\cap K^c\neq \varnothing, $where $D=D(r_{k_1}+ir_{k_2},r_{k_3})$ is the ball with center $r_{k_1}+ir_{k_2} $and ratio $r_{k_3}>0 \}$. Then : \\
$$(E_k)^c=\bigcup\limits_{n \in N,(k_1,k_2,k_3) \in S}L_{n,k_1,k_2,k_3} $$
,where $L_{n,k_1,k_2,k_3}=\{ (g_s)_{s=0}^\infty \in B(K)$, such that there exists a holomorphic function $F$ on $D=D(r_{k_1}+ir_{k_2},r_{k_3})$, such that $F_{\vert D\cap K}={g_k}_{\vert D\cap K}$ and $\vert F(z)\vert\leq n$ for $z \in D$\}. By Baire's Theorem we need to show that $L_{n,k_1,k_2,k_3}$ is closed and nowhere dense. \\
\begin{enumerate}
    \item It is closed, since if we have $g^1,g^2,g^3,.. \in L_{n,k_1,k_2,k_3}$ with $g^n \rightarrow g$ then $g_k^n\rightarrow g_k$ uniformly on K, as $n\rightarrow \infty$, and thus $F_n \rightarrow g_k$ uniformly on $K$, as $n\rightarrow \infty$, where $F_n$ are the extensions of $g^n_k$ on $D$. Since $\vert F_n(z)\vert \leq n$ on $D$, by Montel's Theorem we get a function $F \in H(D)$ such that $F_{s_t}\rightarrow F$ uniformly on compact subsets of $D$ for some sequence $s_t$. Thus $F_{\vert D\cap K}={g_k}_{\vert D\cap K}$ and $\vert F(z)\vert\leq n$ for $z \in D$. Therefore $g \in L_{n,k_1,k_2,k_3}$.
    \item Finally, $L_{n,k_1,k_2,k_3}$ is nowhere dense, since taking $\epsilon>0$, we find $z_0 \in D\cap K^c$ and define $h(z)=\frac{M}{z-z_0}$, where $M$ is small enough such that $d(h,0)<\epsilon$. Then, if $g \in L_{n,k_1,k_2,k_3}$, then $d(g+h,g)<\epsilon$ but $g+h \not\in L_{n,k_1,k_2,k_3}$, because otherwise $h=(h+g)-g \in L_{2n,k_1,k_2,k_3}$, which is not the case. \qedsymbol 
\end{enumerate} 
\section{B(K) in Countable Sets} 

In the last paragraph, we saw various examples of sets of the form $B(K)$, where $K$ is a planar compact set. Notice that in these examples there appears to be a strong correlation between $g_k$ and $g_{k+1}$, where $(g_k)_{k=0}^\infty$ is a sequence in $B(K)$. In this paragraph, we show that this is not always the case. \\
To show our result we will prove a more general theorem, which also holds for non-compact countable sets: \\
\textbf{Theorem 4.1} \textit{Let $X$ denote a countable planar set and $f_0,f_1,f_2,..$ be a sequence of continuous functions on $X$. It is possible to find functions $h_0,h_1,h_2,...$ such that for $i \in N, h_i$ is defined in an open neighborhood of $X, \, h_i$ is locally a polynomial and for $j \in N$ $$\limsup\limits_{i\rightarrow {+\infty},z \in X}\vert h_i^j(z) - f_j(z)\vert \rightarrow 0  $$.} \\
\textbf{Proof.} Let $X=\{x_0,x_1,x_2,x_3,...\}$ and let $k\in N$. For $i=0,1,2,...,k$ we can find an open disk $D_i$ of center $x_i$ and ratio $\delta_i$ such that:
\begin{enumerate}
    \item The disks $D_i$ are disjoint.
    \item There is no point of $X$ lying on the circumference of $D_i$.
    \item $diam(f_j(D_i\cap X)) < \frac{1}{2k}$ for $j=0,1,2,3,...,k$.
    \item If $p_{i,k}$ denotes the unique polynomial of degree $k$ with derivative of order $s$ at $x_i : \,\, f_s(x_i)$ for $s\leq k$, then $diam(p_{i,k}^s(D_i))<\frac{1}{2k}$ for all $i\leq k$ and $s\leq k$.
\end{enumerate} 
Define now the function $h_k$ as follows: 
\\
$h_k=p_{i,k}$ on $D_i$ for every $i\leq k$. In case $X\subset\bigcup\limits_{i=0}^k D_i$, $h_k$ is well defined on an open neighborhood of $X$. If not, pick the least possible $t$ such that $x_t\not\in \bigcup\limits_{i=0}^k D_i$. Then define a disk $D_{k+1}$ of center $x_t$ and ratio $\delta_{k+1}$, disjoint from the already defined disks $D_i, i\leq k$ such that: \\
\begin{enumerate}
    \item There is no point of $X$ lying on the circumference of $D_{k+1}$
    \item $diam(f_j(D_{k+1}\cap X)) < \frac{1}{2k}$ for $j=0,1,2,3,...,k$.
    \item $diam(p_{t,k}^s(D_{k+1}))<\frac{1}{2k}$ for all $s\leq k$.
\end{enumerate} 
Define $h_k=p_{t,k}$ on $D_{k+1}$. In case $X\subset\bigcup\limits_{i=0}^{k+1} D_i$, $h_k$ is well defined on an open neighborhood of $X$. If not, once more pick the least possible $t$ such that $x_t\not\in \bigcup\limits_{i=0}^{k+1} D_i$ and continue this process inductively. \\
We now prove that $h_k^s \rightarrow f_s$ uniformly on $X$ for fixed $s \in N$. Indeed for $k>s$ : $$\sup\limits_{z\in X}\vert h_k^s(z) - f_s(z)\vert=\sup\limits_{i \in N}\vert h_k^s(x_i) - f_s(x_i)\vert $$ 
\begin{enumerate}
    \item For $i\leq k$ we have that $h_k^s(x_i)=f_s(x_i)$
    \item For $i>k \,\,\, x_i$ belongs to some disk $D$ of center $x_p$ having the properties listed above. Therefore : $$\begin{aligned}
\vert h_k^s(x_i) - f_s(x_i)\vert \leq \vert h_k^s(x_i) - h_k^s(x_p) \vert &+ \underbrace{\vert h_k^s(x_p) - f_s(x_p) \vert}_0 + \vert f_s(x_p) - f_s(x_i) \vert \\ &\leq diam(h_k^s(D)) + diam(f_s(D\cap X)) \leq \frac{1}{2k} + \frac{1}{2k}= \frac{1}{k} \end{aligned}$$
\end{enumerate}
\qedsymbol \\
We thus have that: \\ \\
\textbf{Corrolary 4.2} \textit{If K is a countable compact planar set, then $\overline{P}=B(K)=[C(K)]^{\aleph_0}$, where $P$ denotes the set of polynomials and its closure is taken in the topology of $B(K)$. }\\ \\
\textbf{Proof.} The fact that $B(K)=[C(K)]^{\aleph_0}$ is immediate from Theorem 4.1. The part that $\overline{P}=B(K)$ follows from Proposition 2.6 , since $\mathbb C_\infty-K$ is connected.  \qedsymbol

\bibliographystyle{plain}
\bibliography{main}
\noindent
National and Kapodestrian University of Athens \\ University of Athens \\ Department of Mathematics \\ Panepistimiopolis, 157, 84 \\ Athens \\ Greece \\ \\
\noindent
E-mail addresses.\\ Sotiris Armeniakos: sotarmen@outlook.com \\ Giorgos Kotsovolis: gk13@princeton.edu \\ Vassili Nestoridis: vnestor@math.uoa.gr

\end{document}